\documentclass[12pt,leqno]{article}
\usepackage{amssymb}
\usepackage[mathscr]{eucal}
\usepackage{amsmath,amssymb,latexsym,theorem,bbm} %
\usepackage{color,url}
\usepackage{appendix}

\newcommand{\NN}{\mathbb{N}}

\newcommand{\RR}{\mathbb{R}}

\newcommand{\ZZ}{\mathbb{Z}}

\newcommand{\tr}{{\widetilde{r}}}

\newcommand{\cA}{{\mathcal A}}

\newcommand{\cD}{{\mathcal D}}

\newcommand{\cF}{{\mathcal F}}

\newcommand{\dd}{\mathrm{d}}

\newcommand{\EE}{\operatorname{\mathbb{E}}}
\newcommand{\PP}{\operatorname{\mathbb{P}}}
\newcommand{\OO}{\operatorname{O}}
\newcommand{\oo}{\operatorname{o}}

\newcommand{\vare}{\varepsilon}

\renewcommand{\mid}{\,|\,}

\renewcommand{\leq}{\leqslant}
\renewcommand{\geq}{\geqslant}

\newcommand{\distre}{\stackrel{\cD}{=}}

\newcommand{\bbone}{\mathbbm{1}}

\newcommand{\proofend}{\hfill\mbox{$\Box$}}

\numberwithin{equation}{section}

\theoremstyle{change} \theorembodyfont{\em}
\newtheorem{Thm}{Theorem.}[section]
\newtheorem{Lem}[Thm]{Lemma.}
\newtheorem{Pro}[Thm]{Proposition.}

\newtheorem{Def}[Thm]{Definition.}

\theorembodyfont{\rm}
\newtheorem{Rem}[Thm]{Remark.}

\begin{document}

\begin{center}
 {\bfseries\Large Regularly varying non-stationary Galton--Watson \\[2mm]
                   processes with immigration} \\[6mm]
 {\sc\large M\'aty\'as $\text{Barczy}^{*,\diamond}$,
            \ Zsuzsanna $\text{B\H{o}sze}^{**}$,
            \ Gyula $\text{Pap}^{**}$}
\end{center}

\vskip0.2cm

\noindent
 * MTA-SZTE Analysis and Stochastics Research Group,
   Bolyai Institute, University of Szeged,
   Aradi v\'ertan\'uk tere 1, H--6720 Szeged, Hungary.

\noindent
 ** Bolyai Institute, University of Szeged,
    Aradi v\'ertan\'uk tere 1, H-6720 Szeged, Hungary.

\noindent e-mail: barczy@math.u-szeged.hu (M. Barczy),
                  Bosze.Zsuzsanna@stud.u-szeged.hu (Zs.\ B\H{o}sze),
                  papgy@math.u-szeged.hu (G. Pap).

\noindent $\diamond$ Corresponding author.



\renewcommand{\thefootnote}{}
\vspace*{0.2cm}
\footnote{Supported by the Hungarian Croatian Intergovernmental S \& T Cooperation Programme
 for 2017-2018 under Grant No.\ 16-1-2016-0027.
M\'aty\'as Barczy is supported by the J\'anos Bolyai Research Scholarship of the Hungarian
 Academy of Sciences.}

\vspace*{-10mm}

\begin{abstract}
We give sufficient conditions on the initial, offspring and immigration distributions under which
 the distribution of a not necessarily stationary Galton--Watson process with immigration is
 regularly varying at any fixed time.
\end{abstract}

\section{Introduction}
\label{section_intro}

Galton--Watson processes with immigration have been frequently used for modeling the sizes of a
 population over time, so a delicate description of their tail behavior is an important question.
In this paper we focus on regularly varying not necessarily stationary Galton--Watson processes
 with immigration, complementing the results of Basrak et al. \cite{BasKulPal} for the stationary
 case.
By a Galton--Watson process with immigration, we mean a stochastic process \ $(X_n)_{n\geq0}$
 \ given by
 \begin{equation}\label{GWI}
  X_n = \sum_{i=1}^{X_{n-1}} \xi_{n,i} + \vare_n , \qquad n \geq 1 ,
 \end{equation}
 where \ $\bigl\{X_0, \, \xi_{n,i}, \, \vare_n : n, i \geq 1 \bigr\}$ \ are supposed to be
 independent non-negative integer-valued random variables,
 \ $\{\xi_{n,i} : n, i \geq 1\}$ \ and \ $\{\vare_n : n \geq 1\}$ \ are supposed to consist of
 identically distributed random variables, respectively, and \ $\sum_{i=1}^0 := 0$.
\ If \ $\vare_n = 0$, \ $n \geq 1$, \ then we say that \ $(X_n)_{n\geq0}$ \ is a Galton--Watson
 process (without immigration).

Basrak et al.\ \cite{BasKulPal} studied stationary Galton--Watson processes with immigration and
 gave conditions under which the stationary distribution is regularly varying.

In the special case of \ $\PP(\xi_{1,1} = \varrho) = 1$ \ with some non-negative integer
 \ $\varrho$, \ $(X_n)_{n\geq 0}$ \ is nothing else but a first order autoregressive process having
 the form \ $X_n = \varrho X_{n-1} + \vare_n$, \ $n \geq 1$.
\ There is a vast literature on the tail behavior of weighted sums of independent and identically
 distributed regularly varying random variables, especially, of first order autoregressive
 processes with regularly varying noises, see, e.g., Embrechts et al.\
 \cite[Appendix A3.3]{EmbKluMik}.
For instance, in the special case mentioned before, Proposition \ref{GWI_vare} with
 \ $\PP(X_0 = 0) = 1$ \ gives the result of Lemma A3.26 in Embrechts et al. \cite{EmbKluMik}, since
 then \ $X_n = \sum_{i=1}^n \varrho^{n-i} \vare_i$, \ $n \geq 1$.

In Section \ref{section_GWI}, we present conditions on the initial, offspring and immigration
 distributions under which the distribution of a not necessarily stationary Galton--Watson
 process with immigration is regularly varying at any fixed time describing the precise tail
 behavior of the distribution in question as well.
The proofs are delicate applications of Fa\"y et al.\ \cite[Proposition 4.3]{GilGAMikSam} (see
 Proposition \ref{FGAMS}), Robert and Segers \cite[Theorem 3.2]{RobSeg} (see Proposition \ref{RS})
 and Denisov et al.\ \cite[Theorems 1 and 7]{DenFosKor} (see Propositions \ref{DFK0} and
 \ref{DFK}).
We close the paper with four appendicies: in Appendix \ref{section_prel} we recall representations
 of  Galton--Watson process without or with immigration; Appendix \ref{App_moments} is devoted to
 higher moments of Galton--Watson processes; in Appendix \ref{App1} we collect some properties of
 regularly varying functions and distributions used in the paper; and in Appendix \ref{App_rvrs} we
 recall the results of Fa\"y et al.\ \cite[Proposition 4.3]{GilGAMikSam}, Robert and Segers
 \cite[Theorem 3.2]{RobSeg}, Denisov et al.\ \cite[Theorems 1 and 7]{DenFosKor} and some of
 their consequences.

Later on, one may also investigate other tail properties such as intermediate regular
 variation.
Motivated by Bloznelis \cite{Blo}, one may study the asymptotic behavior of the so called local
 probabilities \ $\PP(X_n = \ell)$ \ as \ $\ell \to \infty$ \ for any fixed \ $n \in \NN$.

\section{Tail behavior of Galton--Watson processes with immigration}
\label{section_GWI}

Let \ $\ZZ_+$, \ $\NN$, \ $\RR$, \ $\RR_+$ \ and \ $\RR_{++}$ \ denote the set of non-negative
 integers, positive integers, real numbers, non-negative real numbers and positive real numbers,
 respectively.
For \ $x, y \in \RR$, \ we will use the notations \ $x \land y := \min(x, y)$ \ and
 \ $x \lor y := \max(x, y)$.
\ For functions \ $f : \RR_{++} \to \RR_{++}$ \ and \ $g : \RR_{++} \to \RR_{++}$, \  by the
 notation \ $f(x) \sim g(x)$, \ $f(x) = \oo(g(x))$ \ and \ $f(x) = \OO(g(x))$ \ as
 \ $x \to \infty$, \ we mean that \ $\lim_{x\to\infty} \frac{f(x)}{g(x)} = 1$,
 \ $\lim_{x\to\infty} \frac{f(x)}{g(x)} = 0$ \ and
 \ $\limsup_{x\to\infty} \frac{f(x)}{g(x)} < \infty$, \ respectively.
Every random variable will be defined on a probability space \ $(\Omega, \cA, \PP)$.
\ Equality in distribution of random variables is denoted by \ $\distre$.
\ For notational convenience, let \ $\xi$ \ and \ $\vare$ \ be random variables such that
 \ $\xi \distre \xi_{1,1}$ \ and \ $\vare \distre \vare_1$, and put
 \ $m_\xi := \EE(\xi) \in [0, \infty]$ \ and \ $m_\vare := \EE(\vare) \in [0, \infty]$.

First, we consider the case of regularly varying offspring distribution.

\begin{Pro}\label{GWI_xi}
Let \ $(X_n)_{n\in\ZZ_+}$ \ be a Galton--Watson process with immigration such that \ $\xi$
 \ is regularly varying with index \ $\alpha \in [1, \infty)$ \ and there exists
 \ $r \in (\alpha, \infty)$ \ with \ $\EE(X_0^r) < \infty$ \ and \ $\EE(\vare^r) < \infty$.
\ Suppose that \ $\PP(X_0 = 0) < 1$ \ or \ $\PP(\vare = 0) < 1$.
\ In case of \ $\alpha = 1$, \ assume additionally that \ $m_\xi \in \RR_{++}$.
\ Then for each \ $n \in \NN$, \ we have
 \[
   \PP(X_n > x)
   \sim \EE(X_0) m_\xi^{n-1} \sum_{i=0}^{n-1} m_\xi^{i(\alpha-1)} \PP(\xi > x)
        + m_\vare \sum_{i=1}^{n-1} m_\xi^{n-i-1} \sum_{j=0}^{n-i-1} m_\xi^{j(\alpha-1)} \PP(\xi > x)
 \]
 as \ $x \to \infty$, \ and hence \ $X_n$ \ is also regularly varying with index \ $\alpha$.
\end{Pro}

\noindent{\bf Proof.}
Note that we always have \ $\EE(X_0) \in \RR_+$, \ $m_\xi \in \RR_{++}$ \ and
 \ $m_\vare \in \RR_+$.
\ We use the representation \eqref{GWI_additive}.
Recall that \ $\bigl\{V^{(n)}(X_0), V_i^{(n-i)}(\vare_i) : i \in \{1, \ldots, n\}\bigr\}$
 \ are independent random variables such that \ $V^{(n)}(X_0)$ \ represents the number of
 individuals alive at time \ $n$, \ resulting from the initial individuals \ $X_0$ \ at time \ $0$,
 \ and for each \ $i \in \{1, \ldots, n\}$, \ $V_i^{(n-i)}(\vare_i)$ \ represents the number of
 individuals alive at time \ $n$, \ resulting from the immigration \ $\vare_i$ \ at time \ $i$.
\ If \ $\PP(X_0 = 0) = 1$, \ then \ $\PP(V^{(n)}(X_0) = 0) = 1$, \ otherwise, by Proposition
 \ref{DFK_corr0}, we obtain
 \[
   \PP(V^{(n)}(X_0) > x)
   \sim \EE(X_0) m_\xi^{n-1} \sum_{i=0}^{n-1} m_\xi^{i(\alpha-1)} \PP(\xi > x) \qquad
   \text{as \ $x \to \infty$}
 \]
 once we show
 \begin{equation}\label{GW_VV}
  \PP(V_n > x) \sim m_\xi^{n-1} \sum_{i=0}^{n-1} m_\xi^{i(\alpha-1)} \PP(\xi > x) \qquad
  \text{as \ $x \to \infty$,}
 \end{equation}
 where \ $(V_k)_{k\in\ZZ_+}$ \ is a Galton--Watson process (without immigration) with initial
 value \ $V_0 = 1$ \ and with the same offspring distribution as \ $(X_k)_{k\in\ZZ_+}$.
\ We proceed by induction on \ $n$.
\ For \ $n = 1$, \ \eqref{GW_VV} follows readily, since \ $V_1 = \xi_{1,1}  \distre \xi$.
\ Now let us assume that \eqref{GW_VV} holds for \ $1, \ldots, n - 1$, \ where \ $n \geq 2$.
\ Since \ $(V_k)_{k\in\ZZ_+}$ \ is a time homogeneous Markov process with \ $V_1 = \xi_{1,1}$, \ we
 have \ $V_n \distre V^{(n-1)}(\xi_{1,1})$, \ where \ $(V^{(k)}(\xi_{1,1}))_{k\in\ZZ_+}$ \ is a
 Galton--Watson process (without immigration) with initial value \ $V^{(0)}(\xi_{1,1}) = \xi_{1,1}$
 \ and with the same offspring distribution as \ $(X_k)_{k\in\ZZ_+}$.
\ Applying again the additive property \eqref{GW_additive}, we obtain
 \[
   V_n \distre V^{(n-1)}(\xi_{1,1}) \distre \sum_{i=1}^{\xi_{1,1}} \zeta_i^{(n-1)} ,
 \]
 where \ $\{\zeta_i^{(n-1)} : i \in \NN\}$ \ are independent copies of \ $V_{n-1}$ \ such that
 \ $\{\xi_{1,1}, \zeta_i^{(n-1)} : i \in \NN\}$ \ are independent.
First note that \ $\PP(\zeta_1^{(n-1)} > x) = \OO(\PP(\xi > x))$ \ as \ $x \to \infty$.
\ Indeed, using the induction hypothesis, we obtain
 \[
   \limsup_{x\to\infty} \frac{\PP(\zeta_1^{(n-1)} > x)}{\PP(\xi > x)}
   = \limsup_{x\to\infty} \frac{\PP(V_{n-1} > x)}{\PP(\xi > x)}
   = m_\xi^{n-2} \sum_{i=0}^{n-2} m_\xi^{i(\alpha-1)}
   < \infty .
 \]
Now we can apply Proposition \ref{DFK_corr}, and we obtain
 \begin{align*}
  \PP(V_n > x)
  = \PP\left(\sum_{i=1}^{\xi_{1,1}} \zeta_i^{(n-1)} > x\right)
  &\sim \EE(\xi_{1,1}) \PP(\zeta_1^{(n-1)} > x)
        + (\EE(\zeta_1^{(n-1)}))^\alpha \PP(\xi_{1,1} > x) \\
  &\sim m_\xi \PP(V_{n-1} > x) + m_\xi^{(n-1)\alpha} \PP(\xi > x)
 \end{align*}
 as \ $x \to \infty$, \ since, by \eqref{GWI_EX},
 \ $\EE(\zeta_1^{(n-1)}) = m_\xi^{n-1} \in \RR_{++}$.
\ Using the induction hypothesis and
 \ $m_\xi^{(n-1)\alpha} = m_\xi^{n-1} m_\xi^{(n-1)(\alpha-1)}$, \ we conclude \eqref{GW_VV}.

If \ $\PP(\vare = 0) = 1$, \ then for each \ $i \in \{1, \ldots, n\}$,
 \ $\PP(V_i^{(n-i)}(\vare_i) = 0) = 1$.
\ Otherwise, for each \ $i \in \{1, \ldots, n - 1\}$, \ by \eqref{GW_VV} and Proposition
 \ref{DFK_corr0}, we obtain \ $\PP(V_i^{(n-i)}(\vare_i) > x) \sim m_\vare \PP(V_{n-i} > x)$ \ as
 \ $x \to \infty$, \ and \ $V_n^{(0)}(\vare_n) = \vare_n$.

Applying the convolution property in Lemma \ref{Lem_conv} and \eqref{GW_VV}, we conclude the
 statement.
\proofend

The result of Proposition \ref{GWI_xi} in the special case of \ $\PP(X_0 = 1) = 1$,
 \ $\alpha \in (1, 2)$, \ $m_\xi \in (0, 1)$ \ and \ $\PP(\vare = 0) = 1$ \ has already been
 derived by Basrak et al.\ \cite[page 426]{BasKulPal} written it in an equivalent form
 \[
   \PP(X_n > x) \sim m_\xi^{(n-1)\alpha} \sum_{i=0}^{n-1} m_\xi^{i(1-\alpha)} \PP(\xi > x)
   \qquad \text{as \ $x \to \infty$}
 \]
 for all \ $n \in \NN$.

Note that Denisov et al.\ \cite[Corollary 2]{DenFosKor} proved that
 \[
   \PP(X_2 > x) \sim m_\xi \PP(\xi > x) + \PP\biggl(\xi > \frac{x}{m_\xi}\biggr) \qquad
   \text{as \ $x \to \infty$}
 \]
 if \ $(X_n)_{n\in\ZZ_+}$ \ is a Galton--Watson process (without immigration) such that
 \ $X_0 = 1$, \ $\xi$ \ is intermediate regularly varying \ and \ $m_\xi \in \RR_{++}$, \ and, by
 induction arguments, for each \ $n \in \NN$, \ $\PP(X_n > x) \sim n \PP(\xi > x)$ \ if, in
 addition, \ $m_\xi = 1$.
\ Further, Wachtel et al.\ \cite[formula (5.1)]{WacDenKor} mentioned that for each
 \ $n \in \NN$,
 \begin{equation}\label{WDK(5.1)}
  \PP\biggl(\frac{X_n}{m_\xi^n} > x\biggr) \sim \sum_{i=0}^{n-1} m_\xi^i \PP(\xi > m_\xi^{i+1} x)
  \qquad \text{as \ $x \to \infty$}
 \end{equation}
  if \ $(X_n)_{n\in\ZZ_+}$ \ is a Galton--Watson process (without immigration) such that
 \ $X_0 = 1$, \ $\xi$ \ is intermediate regularly varying \ and \ $m_\xi > 1$.
\ These results for regularly varying \ $\xi$ \ are consequences of Proposition \ref{GWI_xi}.
In fact, Wachtel et al.\ \cite[Theorem 1]{WacDenKor} showed that \eqref{WDK(5.1)} holds uniformly
 in \ $n$, \ and, in particular,
 \[
   \PP\biggl(\frac{X_n}{m_\xi^n} > x\biggr)
   \sim \sum_{i=0}^\infty m_\xi^i \PP(\xi > m_\xi^{i+1} x)
   \qquad \text{as \ $x, n \to \infty$.}
 \]

Next, we consider the case of regularly varying initial distribution.

\begin{Pro}\label{GWI_X_0}
Let \ $(X_n)_{n\in\ZZ_+}$ \ be a Galton--Watson process with immigration such that \ $X_0$
 \ is regularly varying with index \ $\beta \in \RR_+$, \ $\PP(\xi = 0) < 1$ \ and there
 exists \ $r \in (1 \lor \beta, \infty)$ \ with \ $\EE(\xi^r) < \infty$ \ and
 \ $\EE(\vare^r) < \infty$.
\ Then for each \ $n \in \NN$, \ we have
 \[
   \PP(X_n > x) \sim m_\xi^{n\beta} \PP(X_0 > x) \qquad \text{as \ $x \to \infty$,}
 \]
 and hence \ $X_n$ \ is also regularly varying with index \ $\beta$.
\end{Pro}

\noindent{\bf Proof.}
Let us fix \ $n \in \NN$.
\ We use the representation \eqref{GWI_additive}.
In view of the convolution property in Lemma \ref{Lem_conv}, it is enough to prove
 \[
   \PP(V^{(n)}(X_0) > x) \sim m_\xi^{n\beta} \PP(X_0 > x) \qquad \text{as \ $x \to \infty$,}
 \]
 since by Lemma \ref{Lem_seged_momentr}, for each \ $i \in \{1, \ldots, n\}$, \ we have
 \ $\EE((V_i^{(n-i)}(\vare_i))^r) < \infty$ \ yielding
 \ $\EE((\sum_{i=1}^n V_i^{(n-i)}(\vare_i))^r) < \infty$.
\ By the additive property \eqref{GW_additive}, we have
 \ $V^{(n)}(X_0) \distre \sum_{i=1}^{X_0} \zeta_i^{(n)}$.
\ By \eqref{GWI_EX}, \ $\EE(\zeta_1^{(n)}) = m_\xi^n \in \RR_{++}$, \ and by Lemma
 \ref{Lem_seged_momentr}, we have \ $\EE((\zeta_1^{(n)})^r) < \infty$.
\ The statement is a consequence of Proposition \ref{FGAMSRS}.
\proofend

\begin{Pro}\label{GWI_X_0_xi}
Let \ $(X_n)_{n\in\ZZ_+}$ \ be a Galton--Watson process with immigration such that \ $X_0$
 \ and \ $\xi$ \ are regularly varying with index \ $\beta \in [1, \infty)$ \ and
 \ $\PP(\xi > x) = \OO(\PP(X_0 > x))$ \ as \ $x \to \infty$ \ and there exists
 \ $r \in (\beta, \infty)$ \ such that \ $\EE(\vare^r) < \infty$.
\ In case of \ $\beta = 1$, \ assume additionally that \ $\EE(X_0) \in \RR_{++}$ \ and
 \ $m_\xi \in \RR_{++}$.
\ Then for each \ $n \in \NN$, \ we have
 \begin{align*}
  \PP(X_n > x)
   &\sim \EE(X_0) m_\xi^{n-1} \sum_{i=0}^{n-1} m_\xi^{i(\beta-1)} \PP(\xi > x)
         + m_\xi^{n\beta} \PP(X_0 > x) \\
   &\quad
         + m_\vare
           \sum_{i=1}^{n-1} m_\xi^{n-i-1} \sum_{j=0}^{n-i-1} m_\xi^{j(\beta-1)} \PP(\xi > x)
 \end{align*}
 as \ $x \to \infty$, \ and hence \ $X_n$ \ is also regularly varying with index \ $\beta$.
\end{Pro}

\noindent{\bf Proof.}
Let us fix \ $n \in \NN$.
\ Note that we always have \ $\EE(X_0) \in \RR_{++}$ \ and \ $m_\xi \in \RR_{++}$.
\ We use the representation \eqref{GWI_additive}.
By \eqref{GW_VV}, \ $V_n$ \ is regularly varying with index \ $\beta$.
\ By the assumption \ $\PP(\xi > x) = \OO(\PP(X_0 > x))$ \ as \ $x \to \infty$ \ and \eqref{GW_VV},
 we conclude \ $\PP(V_n > x) = \OO(\PP(X_0 > x))$ \ as \ $x \to \infty$.
\ By Proposition \ref{DFK_corr} and \ $\EE(V_n) = m_\xi^n \in \RR_{++}$, \ we obtain
 \[
   \PP(V^{(n)}(X_0) > x) \sim \EE(X_0) \PP(V_n > x) + (\EE(V_n))^\beta \PP(X_0 > x) \qquad
   \text{as \ $x \to \infty$.}
 \]
If \ $\PP(\vare = 0) = 1$, \ then for each \ $i \in \{1, \ldots, n\}$,
 \ $\PP(V_i^{(n-i)}(\vare_i) = 0) = 1$.
\ Otherwise, for each \ $i \in \{1, \ldots, n - 1\}$, \ by \eqref{GW_VV} and Proposition
 \ref{DFK_corr0}, we obtain \ $\PP(V_i^{(n-i)}(\vare_i) > x) \sim m_\vare \PP(V_{n-i} > x)$ \ as
 \ $x \to \infty$, \ and \ $V_n^{(0)}(\vare_n) = \vare_n$.
\ Applying \eqref{GW_VV} and the convolution property, we conclude the statement.
\proofend

Now, we consider the case of regularly varying immigration distribution.

\begin{Pro}\label{GWI_vare}
Let \ $(X_n)_{n\in\ZZ_+}$ \ be a Galton--Watson process with immigration such that \ $\vare$
 \ is regularly varying with index \ $\gamma \in \RR_+$, \ $\PP(\xi = 0) < 1$ \ and there exists
 \ $r \in (1 \lor \gamma, \infty)$ \ with \ $\EE(\xi^r) < \infty$ \ and \ $\EE(X_0^r) < \infty$.
\ Then for each \ $n \in \NN$, \ we have
 \[
   \PP(X_n > x) \sim \sum_{i=1}^n m_\xi^{(n-i)\gamma} \PP(\vare > x) \qquad
   \text{as \ $x \to \infty$,}
 \]
 and hence \ $X_n$ \ is also regularly varying with index \ $\gamma$.
\end{Pro}

\noindent{\bf Proof.}
Let us fix \ $n \in \NN$.
\ We use the representation \eqref{GWI_additive}.
By Lemma \ref{Lem_seged_momentr}, we have \ $\EE((V^{(n)}(X_0))^r) < \infty$.
\ For each \ $i \in \{1, \ldots, n\}$, \ applying Proposition \ref{GWI_X_0} with \ initial
 distribution \ $\vare$ \ and without immigration, we obtain
 \[
   \PP(V_i^{(n-i)}(\vare_i) > x) \sim m_\xi^{(n-i)\gamma} \PP(\vare > x) \qquad
   \text{as \ $x \to \infty$.}
 \]
By the convolution property in Lemma \ref{Lem_conv}, we conclude the statement.
\proofend

\begin{Pro}\label{GWI_xi_vare}
Let \ $(X_n)_{n\in\ZZ_+}$ \ be a Galton--Watson process with immigration such that \ $\xi$
 \ and \ $\vare$ \ are regularly varying with index \ $\gamma \in [1, \infty)$,
 \ $\PP(\xi > x) = \OO(\PP(\vare > x))$ \ as \ $x \to \infty$ \ and there exists
 \ $r \in (\gamma, \infty)$ \ with \ $\EE(X_0^r) < \infty$.
\ In case of \ $\gamma = 1$, \ assume additionally that \ $m_\xi \in \RR_{++}$ \ and
 \ $m_\vare \in \RR_{++}$.
\ Then for each \ $n \in \NN$, \ we have
\begin{align*}
 \PP(X_n > x)
 &\sim \EE(X_0) m_\xi^{n-1} \sum_{i=0}^{n-1} m_\xi^{i(\gamma-1)} \PP(\xi > x) \\
 &\quad
       + m_\vare \sum_{j=1}^{n-1} m_\xi^{n-j-1} \sum_{i=0}^{n-j-1} m_\xi^{i(\gamma-1)} \PP(\xi > x)
       + \sum_{j=1}^n m_\xi^{(n-j)\gamma} \PP(\vare > x)
\end{align*}
 as \ $x \to \infty$, \ and hence \ $X_n$ \ is also regularly varying with index \ $\gamma$.
\end{Pro}

\noindent{\bf Proof.}
Let us fix \ $n \in \NN$.
\ Note that we always have \ $\EE(X_0) \in \RR_+$, \ $m_\xi \in \RR_{++}$ \ and
 \ $m_\vare \in \RR_{++}$.
\ We use the representation \eqref{GWI_additive}.
If \ $\PP(X_0 = 0) = 1$, \ then \ $\PP(V^{(n)}(X_0) = 0) = 1$, \ otherwise, applying
 Proposition \ref{GWI_xi} without immigration, we obtain
 \[
   \PP(V^{(n)}(X_0) > x)
   \sim \EE(X_0) m_\xi^{n-1} \sum_{i=0}^{n-1} m_\xi^{i(\gamma-1)} \PP(\xi > x) \qquad
   \text{as \ $x \to \infty$.}
 \]
For each \ $i \in \{1, 2, \ldots, n - 1\}$, \ by \eqref{GW_VV}, \ $V_{n-i}$ \ is regularly varying with
 index \ $\gamma$, \ and hence, by Proposition \ref{DFK_corr}, we get
 \[
   \PP(V^{(n-i)}_i(\vare_i) > x)
   \sim m_\vare \PP(V_{n-i} > x) + (\EE(V_{n-i}))^\gamma \PP(\vare > x) \qquad
   \text{as \ $x \to \infty$.}
 \]
Since \ $V^{(0)}_n(\vare_n) = \vare_n$ \ and \ $V_0 = 1$, \ the above asymptotics is valid
 for \ $i = n$ \ as well.
Applying again \eqref{GW_VV} and the convolution property in Lemma \ref{Lem_conv} together with the
 fact that \ $\EE(V_{n-i}) = m_\xi^{n-i}$ \ for all \ $i \in \{1, 2, \ldots, n\}$, \ we conclude
 the statement.
\proofend

\begin{Pro}\label{GWI_X_0_vare}
Let \ $(X_n)_{n\in\ZZ_+}$ \ be a Galton--Watson process with immigration such that \ $X_0$
 \ and \ $\vare$ \ are regularly varying with index \ $\gamma \in \RR_+$, \ $\PP(\xi = 0) < 1$
 \ and there exists \ $r \in (1 \lor \gamma, \infty)$ \ with \ $\EE(\xi^r) < \infty$.
\ Then for each \ $n \in \NN$, \ we have
 \[
   \PP(X_n > x) \sim m_\xi^{n\gamma} \PP(X_0 > x) + \sum_{i=1}^n m_\xi^{(n-i)\gamma} \PP(\vare > x)
   \qquad \text{as \ $x \to \infty$,}
 \]
 and hence \ $X_n$ \ is also regularly varying with index \ $\gamma$.
\end{Pro}

\noindent{\bf Proof.}
We use the representation \eqref{GWI_additive}.
Applying Proposition \ref{GWI_X_0} first with \ $X_0$ \ and then with \ $X_0 \distre \vare$ \ (in
 both cases without immigration), then the convolution property in Lemma \ref{Lem_conv}, we
 conclude the statement.
\proofend

\begin{Pro}\label{GWI_X_0_xi_vare_3}
Let \ $(X_n)_{n\in\ZZ_+}$ \ be a Galton--Watson process with immigration such that \ $X_0$, \ $\xi$
 \ and \ $\vare$ \ are regularly varying with index \ $\beta \in [1, \infty)$,
 \ $\PP(\xi > x) = \OO(\PP(X_0 > x))$ \ as \ $x \to \infty$ \ and
 \ $\PP(\xi > x) = \OO(\PP(\vare > x))$ \ as \ $x \to \infty$.
\ In case of \ $\beta = 1$, \ assume additionally that \ $\EE(X_0) \in \RR_{++}$,
 \ $m_\xi \in \RR_{++}$ \ and \ $m_\vare \in \RR_+$.
\ Then for each \ $n \in \NN$, \ we have
 \begin{align*}
  \PP(X_n > x)
  &\sim \EE(X_0) m_\xi^{n-1} \sum_{i=0}^{n-1} m_\xi^{i(\beta-1)} \PP(\xi > x)
        + m_\xi^{n\beta} \PP(X_0 > x) \\
  &\quad
        + m_\vare \sum_{j=1}^{n-1} m_\xi^{n-j-1} \sum_{i=0}^{n-j-1} m_\xi^{i(\beta-1)} \PP(\xi > x)
        + \sum_{j=1}^n m_\xi^{(n-j)\beta} \PP(\vare > x)
 \end{align*}
 as \ $x \to \infty$, \ and hence \ $X_n$ \ is also regularly varying with index \ $\beta$.
\end{Pro}

\noindent{\bf Proof.}
Let us fix \ $n \in \NN$.
\ We use the representation \eqref{GWI_additive}.
Applying Proposition \ref{GWI_X_0_xi} first with \ $X_0$ \ and then with \ $X_0 \distre \vare$
 \ (in both cases without immigration), and then the convolution property in Lemma \ref{Lem_conv},
 we conclude the statement.
\proofend

\begin{Rem}
Note that the situation when \ $(X_n)_{n\in\ZZ_+}$ \ is a Galton--Watson process with
 immigration such that \ $\xi$ \ is regularly varying with index \ $\alpha \in [1, \infty)$,
 \ $X_0$ \ is regularly varying with index \ $\beta \in (\alpha, \infty)$,
 \ $\PP(\xi = 0) < 1$ \ and there exists an \ $r \in (\beta, \infty)$ \ such that
 \ $\EE(\vare^r) < \infty$ \ is covered by Proposition \ref{GWI_xi}, since then
 \ $\EE(X_0^\tr) < \infty$ \ for all \ $\tr \in (\alpha, \beta)$.
\ Moreover, the situation when \ $(X_n)_{n\in\ZZ_+}$ \ is a Galton--Watson process with
 immigration such that \ $X_0$ \ is regularly varying with index \ $\beta \in \RR_+$, \ $\xi$ \ is
 regularly varying with index \ $\alpha \in (1 \lor \beta, \infty)$ \ and there exists
 \ $r \in (1 \lor \beta, \infty)$ \ with \ $\EE(\vare^r) < \infty$ \ is covered by Proposition
 \ref{GWI_X_0}, since then \ $\PP(\xi = 0) < 1$ \ and \ $\EE(\xi^\tr) < \infty$ \ for all
 \ $\tr \in (1 \lor \beta, \alpha)$.
\ The case of \ $\alpha = \beta$ \ is considered in Proposition \ref{GWI_X_0_xi}.
One could formulate other special cases of our results.
\proofend
\end{Rem}

\appendix

\vspace*{5mm}

\noindent{\bf\Large Appendices}

\section{Representations of Galton--Watson processes without or with immigration}
\label{section_prel}

If \ $(X_n)_{n\in\ZZ_+}$ \ is a Galton--Watson process (without immigration), then for each
 \ $n \in \NN$, \ the additive (or branching) property of a Galton--Watson process (without
 immigration), see, e.g.\ in Athreya and Ney \cite[Chapter I, Part A, Section 1]{AthNey}, together
 with the law of total probability, imply
 \begin{equation}\label{GW_additive}
  X_n \distre \sum_{i=1}^{X_0} \zeta_i^{(n)} ,
 \end{equation}
 where \ $\{\zeta_i^{(n)} : i \in \NN\}$ \ are independent copies of \ $V_n$ \ such that
 \ $\{X_0, \zeta_i^{(n)} : i \in \NN\}$ \ are independent, and \ $(V_k)_{k\in\ZZ_+}$ \ is a
 Galton--Watson process (without immigration) with initial value \ $V_0 = 1$ \ and with the same
 offspring distribution as \ $(X_k)_{k\in\ZZ_+}$.

If \ $(X_n)_{n\in\ZZ_+}$ \ is a Galton--Watson process with immigration, then for each
 \ $n \in \NN$, \ we have
 \begin{equation}\label{GWI_additive}
  X_n = V^{(n)}(X_0) + \sum_{i=1}^n V_i^{(n-i)}(\vare_i) ,
 \end{equation}
 where \ $\bigl\{V^{(n)}(X_0), V_i^{(n-i)}(\vare_i) : i \in \{1, \ldots, n\}\bigr\}$ \ are
 independent random variables such that \ $V^{(n)}(X_0)$ \ represents the number of individuals
 alive at time \ $n$, \ resulting from the initial individuals \ $X_0$ \ at time \ $0$, \ and for
 each \ $i \in \{1, \ldots, n\}$, \ $V_i^{(n-i)}(\vare_i)$ \ represents the number of individuals
 alive at time \ $n$, \ resulting from the immigration \ $\vare_i$ \ at time \ $i$, \ see, e.g.,
 Kaplan \cite[formula (1.1)]{Kap2}.
Clearly, \ $(V^{(k)}(X_0))_{k\in\ZZ_+}$ \ and \ $(V_i^{(k)}(\vare_i))_{k\in\ZZ_+}$,
 \ $i \in \{1, \ldots, n\}$, \ are independent Galton--Watson processes (without immigration) with
 initial values \ $V^{(0)}(X_0) = X_0$ \ and \ $V_i^{(0)}(\vare_i) = \vare_i$,
 \ $i \in \{1, \ldots, n\}$, \ respectively, with the same offspring distributions as
 \ $(X_k)_{k\in\ZZ_+}$.

\section{Moment estimation for Galton--Watson processes}
\label{App_moments}

Next, we recall some results for the expectation of a Galton--Watson process (without immigration).
\ If \ $m_\xi \in \RR_+$ \ and \ $\EE(X_0) \in \RR_+$, \ then \eqref{GWI} implies
 \ $\EE(X_n \mid \cF_{n-1}) = X_{n-1} m_\xi$, \ $n \in \NN$, \ where
 \ $\cF_n := \sigma(X_0, \dots, X_n)$, \ $n \in \ZZ_+$.
\ Consequently, \ $\EE(X_n) = m_\xi \EE(X_{n-1})$, \ $n \in \NN$, \ thus
 \begin{equation}\label{GWI_EX}
  \EE(X_n) = m_\xi^n \EE(X_0) , \qquad n \in \NN .
 \end{equation}
Next, we present an auxiliary lemma on higher moments of \ $(X_n)_{n\in\ZZ_+}$.

\begin{Lem}\label{Lem_seged_momentr}
Let \ $(X_n)_{n\in\ZZ_+}$ \ be a Galton--Watson process (without immigration) such that
 \ $\EE(X_0^r) < \infty$ \ and \ $\EE(\xi^r) < \infty$ \ with some \ $r > 1$.
\ Then \ $\EE(X_n^r) < \infty$ \ for all \ $n \in \NN$.
\end{Lem}

\noindent{\bf Proof.}
By power means inequality, we have
 \[
   \EE(X_n^r \mid \cF_{n-1})
   = \EE\Biggl(\Biggl(\sum_{i=1}^{X_{n-1}} \xi_{n,i}\Biggr)^r \,\Bigg|\, \cF_{n-1}\Biggr)
   \leq \EE\Biggl(X_{n-1}^{r-1} \sum_{i=1}^{X_{n-1}} \xi_{n,i}^r \,\Bigg|\, \cF_{n-1}\Biggr)
   = X_{n-1}^r \EE(\xi^r)
    < \infty
 \]
 for all \ $n \in \NN$.
\ Consequently, \ $\EE(X_n^r) \leq \EE(\xi^r)\EE(X_{n-1}^r)$, $n\in\NN$, \ thus
 \ $\EE(X_n^r) \leq \EE(X_0^r) (\EE(\xi^r))^n$, \ $n \in \NN$, \ yielding the assertion.
\proofend

\section{Regularly varying distributions}
\label{App1}

First, we recall the notions of slowly varying and regularly varying functions, respectively.

\begin{Def}
A measurable function \ $U: \RR_{++} \to \RR_{++}$ \ is called regularly varying at infinity with
 index \ $\rho \in \RR$ \ if for all \ $q \in \RR_{++}$,
 \[
   \lim_{x\to\infty} \frac{U(qx)}{U(x)} = q^\rho .
 \]
In case of \ $\rho = 0$, \ we call \ $U$ \ slowly varying at infinity.
\end{Def}

Next, we recall the notion of regularly varying non-negative random variables.

\begin{Def}
A non-negative random variable \ $X$ \ is called regularly varying with index \ $\alpha \in \RR_+$
 \ if \ $U(x) := \PP(X > x) \in \RR_{++}$ \ for all \ $x \in \RR_{++}$, \ and \ $U$ \ is regularly
 varying at infinity with index \ $-\alpha$.
\end{Def}

\begin{Def}
Let \ $X$ \ be a non-negative random variable such that \ $\PP(X > x) \in \RR_{++}$ \ for all
 \ $x \in \RR_{++}$.
\ We call \ $X$
 \begin{itemize}
  \item
   long-tailed if \ $\PP(X>x+y) \sim \PP(X>x)$ \ as \ $x \to \infty$ \ for any fixed
    \ $y \in \RR_{++}$;
  \item
   dominated varying if \ $\PP(X>xy) = \OO(\PP(X>x))$ \ as \ $x \to \infty$ \ for all (or, equivalently,
   for some) \ $y \in (0, 1)$;
  \item
   intermediate regularly varying (also called consistently varying) if
    \[
      \lim_{\vare\downarrow0} \limsup_{x\to\infty} \frac{\PP(X>(1-\vare)x)}{\PP(X>x)} = 1;
    \]
  \item
   strongly subexponential if \ $\EE(X) < \infty$ \ and
    \[
      \int_0^x \PP(X > x - y) \PP(X > y) \, \dd y \sim 2 \EE(X) \PP(X > x) \qquad
      \text{as \ $x \to \infty$.}
    \]
 \end{itemize}
\end{Def}

Note that if \ $X$ \ is a non-negative regularly varying random variable with index
 \ $\alpha \in \RR_+$, \ then \ $X$ \ is intermediate regularly varying as well.

\begin{Lem}\label{sv}
If \ $L : \RR_{++} \to \RR_{++}$ \ is a slowly varying function (at infinity), then
 \[
   \lim_{x\to\infty} x^\delta L(x) = \infty ,\qquad
   \lim_{x\to\infty} x^{-\delta} L(x) = 0 , \qquad \delta \in \RR_{++} .
 \]
\end{Lem}
For Lemma \ref{sv}, see, Bingham et al.\ \cite[Proposition 1.3.6. (v)]{BinGolTeu}.

\begin{Lem}\label{rvexp}
If \ $X$ \ is a non-negative regularly varying random variable with index \ $\alpha \in \RR_{++}$,
 \ then \ $\EE(X^\beta) < \infty$ \ for all \ $\beta \in (-\infty, \alpha)$ \ and
 \ $\EE(X^\beta) = \infty$ \ for all \ $\beta \in (\alpha, \infty)$.
\end{Lem}
For Lemma \ref{rvexp}, see, e.g., Embrechts et al.\ \cite[Proposition A3.8]{EmbKluMik}.

\begin{Lem}\label{exposv}
If \ $X$ \ and \ $Y$ \ are non-negative random variables such that \ $X$ \ is regularly varying
 with index \ $\alpha \in \RR_+$ \ and there exists \ $r \in (\alpha, \infty)$ \ with
 \ $\EE(Y^r) < \infty$, \ then \ $\PP(Y > x) = \oo(\PP(X > x))$ \ as \ $x \to \infty$.
\end{Lem}

\noindent{\bf Proof.}
Applying Lemma \ref{sv}, we obtain
 \[
   0 \leq \frac{\PP(Y > x)}{\PP(X > x)}
     \leq \frac{\EE(Y^r)}{x^r\PP(X > x)}
     = \frac{\EE(Y^r)}{x^{r-\alpha}L(x)}
     \to 0 \qquad \text{as \ $x \to \infty$,}
 \]
 where \ $L(x) := x^\alpha \PP(X > x)$, \ $x \in \RR_{++}$, \ is a slowly varying function.
\proofend

Combining Lemmas \ref{rvexp} and \ref{exposv}, we obtain the following corollary.

\begin{Lem}\label{svosv}
If \ $X_1$ \ and \ $X_2$ \ are non-negative regularly varying random variables with index
 \ $\alpha_1 \in \RR_+$ \ and \ $\alpha_2 \in \RR_+$, \ respectively, such that
 \ $\alpha_1 < \alpha_2$, \ then \ $\PP(X_2 > x) = \oo(\PP(X_1 > x))$ \ as \ $x \to \infty$.
\end{Lem}

\begin{Lem}[Convolution property]\label{Lem_conv}
If \ $X_1$ \ and \ $X_2$ \ are non-negative random variables such that \ $X_1$ \ is regularly
 varying with index \ $\alpha \in \RR_+$ \ and \ there exists \ $r \in (\alpha, \infty)$ \ with
 \ $\EE(X_2^r) < \infty$, \ then \ $\PP(X_1 + X_2 > x) \sim \PP(X_1 > x)$ \ as \ $x \to \infty$,
 \ and hence \ $X_1 + X_2$ \ is regularly varying with index \ $\alpha$.

If \ $X_1$ \ and \ $X_2$ \ are independent non-negative regularly varying random variables with
 index \ $\alpha \in \RR_+$, \ then \ $\PP(X_1 + X_2 > x) \sim \PP(X_1 > x) + \PP(X_2 > x)$ \ as
 \ $x \to \infty$, \ and hence \ $X_1 + X_2$ \ is regularly varying with index \ $\alpha$.
\end{Lem}
The statements of Lemma \ref{Lem_conv} follow, e.g., from parts 1 and 3 of Lemma B.6.1 of
 Buraczewski et al.\ \cite{BurDamMik} and Lemmas \ref{exposv} and \ref{svosv} together with the
 fact that the sum of two slowly varying functions is slowly varying.

\begin{Thm}[Karamata's theorem for truncated moments]\label{Krtrthm}
Consider a non-negati\-ve regularly varying random variable \ $X$ \ with index
 \ $\alpha \in \RR_{++}$.
\ Then we have
 \begin{align*}
  &\lim_{x\to\infty} \frac{x^\beta\PP(X>x)}{\EE(X^\beta\bbone_{\{X\leq x\}})}
   = \frac{\beta-\alpha}{\alpha} \qquad \text{for \ $\beta \in (\alpha, \infty)$,} \\
  &\lim_{x\to\infty} \frac{x^\beta\PP(X>x)}{\EE(X^\beta\bbone_{\{X>x\}})}
   = \frac{\alpha-\beta}{\alpha} \qquad \text{for \ $\beta \in (-\infty, \alpha)$.}
 \end{align*}
\end{Thm}

For Theorem \ref{Krtrthm}, see, e.g., Bingham et al.\ \cite[pages 26--27]{BinGolTeu} or
 Buraczewski et al.\ \cite[Appendix B.4]{BurDamMik}.

\section{Regularly varying random sums}
\label{App_rvrs}

Now, we recall sufficient conditions under which a random sum is regularly varying.

\begin{Pro}\label{FGAMS}
Let \ $\tau$ \ be a non-negative integer-valued random variable and let
 \ $\{\zeta, \zeta_i : i \in \NN\}$ \ be independent and identically distributed non-negative
 random variables, independent of \ $\tau$, \ such that \ $\tau$ \ is regularly varying with index
 \ $\beta \in [0, 1)$ \ and \ $\EE(\zeta) \in \RR_{++}$.
\ Then we have
 \[
   \PP\biggl(\sum_{i=1}^\tau \zeta_i > x\biggr)
   \sim \PP\biggl(\tau > \frac{x}{\EE(\zeta)}\biggr)
   \sim (\EE(\zeta))^\beta \PP(\tau > x) \qquad \text{as \ $x \to \infty$.}
 \]
\end{Pro}

Proposition \ref{FGAMS} follows from Proposition 4.3 in Fa\"y et al.\ \cite{GilGAMikSam},
 since in case of \ $\beta \in [0, 1)$, \ the condition \ $\PP(\zeta > x) = \oo(\PP(\tau > x))$
 \ as \ $x \to \infty$ \ is automatically satisfied, see Lemma \ref{exposv}.
Fa\"y et al.\ \cite[Proposition 4.3]{GilGAMikSam} claim the same result for \ $\beta = 1$ \ under
 the additional assumption \ $\EE(\tau) < \infty$ \ and also for \ $\beta \in (1, \infty)$, \ but
 their proof is not complete, so for \ $\beta \in [1, \infty)$ \ we will use the following result
 of Robert and Segers \cite[Theorem 3.2]{RobSeg}.

\begin{Pro}\label{RS}
Let \ $\tau$ \ be a non-negative integer-valued random variable and let
 \ $\{\zeta, \zeta_i : i \in \NN\}$ \ be independent and identically distributed non-negative
 random variables, independent of \ $\tau$, \ such that \ $\tau$ \ is intermediate regularly
 varying, \ $\EE(\zeta) \in \RR_{++}$ \ and there exists \ $r \in (1, \infty)$ \ with
 \ $\EE(\zeta^r) < \infty$.
\ Assume that one of the following two conditions holds:
 \begin{itemize}
  \item
   $\EE(\tau) < \infty$ \ and \ $\PP(\zeta > x) = \oo(\PP(\tau > x))$ \ as \ $x \to \infty$;
  \item
   $\EE(\tau) = \infty$ \ and there exists \ $q\in[1,r)$ \ such that
   \ $\limsup_{x\to\infty} \frac{\EE(\tau\bbone_{\{\tau\leq x\}})}{x^q\PP(\tau>x)} < \infty$.
 \end{itemize}
Then we have
 \[
   \PP\biggl(\sum_{i=1}^\tau \zeta_i > x\biggr)
   \sim \PP\biggl(\tau > \frac{x}{\EE(\zeta)}\biggr) \qquad \text{as \ $x \to \infty$.}
 \]
\end{Pro}

Combining Propositions \ref{FGAMS} and \ref{RS}, we obtain the following result.

\begin{Pro}\label{FGAMSRS}
Let \ $\tau$ \ be a non-negative integer-valued random variable and let
 \ $\{\zeta, \zeta_i : i \in \NN\}$ \ be independent and identically distributed non-negative
 random variables, independent of \ $\tau$, \ such that \ $\tau$ \ is regularly varying with index
 \ $\beta \in \RR_+$ \ and \ $\EE(\zeta) \in \RR_{++}$.
\ In case of \ $\beta \in [1,\infty)$, \ assume additionally that there exists
 \ $r \in (\beta, \infty)$ \ with \ $\EE(\zeta^r) < \infty$.
\ Then we have
 \[
   \PP\biggl(\sum_{i=1}^\tau \zeta_i > x\biggr)
   \sim \PP\biggl(\tau > \frac{x}{\EE(\zeta)}\biggr)
   \sim (\EE(\zeta))^\beta \PP(\tau > x) \qquad \text{as \ $x \to \infty$,}
 \]
 and hence \ $\sum_{i=1}^\tau \zeta_i$ \ is also regularly varying with index \ $\beta$.
\end{Pro}

\noindent{\bf Proof.}
By Lemma \ref{exposv}, we have \ $\PP(\zeta > x) = \oo(\PP(\tau > x))$ \ as
 \ $x \to \infty$.

In case of \ $\beta \in [0, 1)$, \ the statement follows by Proposition \ref{FGAMS}.

In case of \ $\beta \in [1, \infty)$ \ and \ $\EE(\tau) < \infty$, \ the statement follows by
 Proposition \ref{RS}.
Indeed, any regularly varying random variable is intermediate regularly varying, hence
 \ $\tau$ \ is intermediate regularly varying.

In case of \ $\beta \in [1, \infty)$ \ and \ $\EE(\tau) = \infty$, \ we have \ $\beta = 1$, \ and
 by Karamata's theorem for truncated moments (see Theorem \ref{Krtrthm}), for each
 \ $q \in (1, \infty)$,
 \[
   \lim_{x\to\infty} \frac{x^q\PP(\tau>x)}{\EE(\tau^q\bbone_{\{\tau\leq x\}})} = q - 1 ,
 \]
 hence
 \[
   \frac{\EE(\tau\bbone_{\{\tau\leq x\}})}{x^q\PP(\tau>x)}
   = \frac{\EE(\tau^q\bbone_{\{\tau\leq x\}})}{x^q\PP(\tau>x)}
     \frac{\EE(\tau\bbone_{\{\tau\leq x\}})}{\EE(\tau^q\bbone_{\{\tau\leq x\}})}
   \to \frac{1}{q-1} \cdot 0
   = 0 \qquad \text{as \ $x \to \infty$,}
 \]
 since
 \[
   0 \leq \frac{\EE(\tau\bbone_{\{\tau\leq x\}})}{\EE(\tau^q\bbone_{\{\tau\leq x\}})}
     \leq \frac{\EE(\tau\bbone_{\{\tau\leq x\}})}
               {\EE(\tau^q\bbone_{\{\frac{x}{2}\leq\tau\leq x\}})}
     \leq \frac{x}{(\frac{x}{2})^q}
     = 2^q x^{1-q}
     \to 0 \qquad \text{as \ $x \to \infty$,}
 \]
 which can be also found in remark after Theorem 3.2 in Robert and Segers \cite{RobSeg}.
Since \ $r > \beta \geq 1$, \ by Proposition \ref{RS}, we have the statement.
\proofend

The next proposition is a special case of part (ii) of Theorem 1 in Denisov et al.\
 \cite{DenFosKor}.

\begin{Pro}\label{DFK0}
Let \ $\tau$ \ be a non-negative integer-valued random variable and let
 \ $\{\zeta, \zeta_i : i \in \NN\}$ \ be independent and identically distributed non-negative
 random variables, independent of \ $\tau$, \ such that \ $\zeta$ \ is strongly subexponential
 (yielding $\EE(\zeta) \in \RR_{++}$), \ $\EE(\tau) \in \RR_{++}$ \ and there exists
 \ $c \in (\EE(\zeta), \infty)$ \ with \ $\PP(c \tau > x) = \oo(\PP(\zeta > x))$ \ as
 \ $x \to \infty$.
\ Then we have
 \[
   \PP\biggl(\sum_{i=1}^\tau \zeta_i > x\biggr) \sim \EE(\tau) \PP(\zeta > x) \qquad
   \text{as \ $x \to \infty$.}
 \]
\end{Pro}

As a consequence, we obtain the following corollary.

\begin{Pro}\label{DFK_corr0}
Let \ $\tau$ \ be a non-negative integer-valued random variable and let
 \ $\{\zeta, \zeta_i : i \in \NN\}$ \ be independent and identically distributed non-negative
 random variables, independent of \ $\tau$, \ such that \ $\zeta$ \ is regularly varying with index
 \ $\alpha \in [1, \infty)$, \ $\PP(\tau = 0) < 1$ \ and there exists \ $r \in (\alpha, \infty)$
 \ with \ $\EE(\tau^r) < \infty$.
\ In case of \ $\alpha = 1$, \ assume additionally that \ $\EE(\zeta) \in \RR_{++}$.
\ Then we have
 \[
   \PP\biggl(\sum_{i=1}^\tau \zeta_i > x\biggr) \sim \EE(\tau) \PP(\zeta > x) \qquad
   \text{as \ $x \to \infty$,}
 \]
 and hence \ $\sum_{i=1}^\tau \zeta_i$ \ is also regularly varying with index \ $\alpha$.
\end{Pro}

\noindent{\bf Proof.}
Note that we always have \ $\EE(\tau) \in \RR_{++}$ \ and \ $\EE(\zeta) \in \RR_{++}$.
\ Any regularly varying random variable is intermediate regularly varying, hence \ $\zeta$ \ is
 intermediate regularly varying.
Any intermediate regularly varying distribution is long-tailed and dominated varying, thus
 \ $\zeta$ \ is long-tailed and dominated varying.
Taking into account that \ $\EE(\zeta) < \infty$, \ $\zeta$ \ is strongly subexponential, see
 Kl\"uppelberg \cite[Theorem 3.2 (a)]{Klu}.
By Lemma \ref{exposv}, \ $\EE(\tau^r) < \infty$ \ implies \ $\PP(c \tau > x) = \oo(\PP(\zeta > x))$
 \ as \ $x \to \infty$ \ for any \ $c \in (\EE(\zeta), \infty)$, \ hence Proposition \ref{DFK0}
 yields the statement.
\proofend

Note that the situation when \ $\tau$ \ is regularly varying with index \ $\beta \in [1, \infty)$,
 \ $\zeta$ \ is regularly varying with index \ $\alpha \in (\beta, \infty)$ \ and
 \ $\EE(\tau) \in \RR_{++}$ \ is covered by Proposition
 \ref{FGAMSRS}, since then \ $\EE(\zeta^r) < \infty$ \ for all \ $r \in (\beta, \alpha)$.
\ Moreover, the situation when \ $\tau$ \ is regularly varying with index \ $\beta \in (1, \infty)$,
 \ $\zeta$ \ is regularly varying with index \ $\alpha \in [1, \beta)$ \ and \ $\EE(\zeta) \in \RR_{++}$
 \ is covered by Proposition \ref{DFK_corr0}, since then \ $\EE(\tau^r) < \infty$ \ for all
 \ $r \in (\alpha, \beta)$.
The case \ $\alpha = \beta$ \ will be covered by a corollary of the next proposition, which is due
 to Denisov et al.\ \cite[Theorem 7]{DenFosKor}.

\begin{Pro}\label{DFK}
Let \ $\tau$ \ be a non-negative integer-valued random variable and let
 \ $\{\zeta, \zeta_i : i \in \NN\}$ \ be independent and identically distributed non-negative
 random variables, independent of \ $\tau$, \ such that \ $\tau$ \ is intermediate regularly
 varying, \ $\EE(\tau) \in \RR_{++}$, \ $\PP(\zeta > x) = \OO(\PP(\tau > x))$ \ as \ $x \to \infty$
 \ and \ $\zeta$ \ is strongly subexponential (yielding $\EE(\zeta) \in \RR_{++}$).
\ Then we have
 \[
   \PP\biggl(\sum_{i=1}^\tau \zeta_i > x\biggr)
   \sim \EE(\tau) \PP(\zeta > x) + \PP\biggl(\tau > \frac{x}{\EE(\zeta)}\biggr) \qquad
   \text{as \ $x \to \infty$.}
 \]
\end{Pro}

As a consequence, we obtain the following corollary.

\begin{Pro}\label{DFK_corr}
Let \ $\tau$ \ be a non-negative integer-valued random variable and let
 \ $\{\zeta, \zeta_i : i \in \NN\}$ \ be independent and identically distributed non-negative
 random variables, independent of \ $\tau$, \ such that \ $\tau$ \ and \ $\zeta$ \ are regularly
 varying with index \ $\beta \in [1, \infty)$, \ and \ $\PP(\zeta > x) = \OO(\PP(\tau > x))$ \ as
 \ $x \to \infty$.
\ In case of \ $\beta = 1$, \ assume additionally that \ $\EE(\tau) \in \RR_{++}$ \ and
 \ $\EE(\zeta) \in \RR_{++}$.
\ Then we have
 \[
   \PP\biggl(\sum_{i=1}^\tau \zeta_i > x\biggr)
   \sim \EE(\tau) \PP(\zeta > x) + (\EE(\zeta))^\beta \PP(\tau > x) \qquad
   \text{as \ $x \to \infty$,}
 \]
 and hence \ $\sum_{i=1}^\tau \zeta_i$ \ is also regularly varying with index \ $\beta$.
\end{Pro}

\noindent{\bf Proof.}
Note that we always have \ $\EE(\tau) \in \RR_{++}$ \ and \ $\EE(\zeta) \in \RR_{++}$.
\ Any regularly varying random variable is intermediate regularly varying, hence \ $\tau$ \ and
 \ $\zeta$ \ are intermediate regularly varying.
Any intermediate regularly varying distribution is long-tailed and dominated varying, thus
 \ $\zeta$ \ is long-tailed and dominated varying.
Taking into account that \ $\EE(\zeta) < \infty$, \ $\zeta$ \ is strongly subexponential, see
 Kl\"uppelberg \cite[Theorem 3.2 (a)]{Klu}, hence Proposition \ref{DFK} yields the statement.
\proofend

Note that in the situation when \ $\tau$ \ is regularly varying with index
 \ $\beta \in [1, \infty)$, \ $\zeta$ \ is regularly varying with index \ $\alpha \in [1, \beta)$,
 \ $\EE(\tau) \in \RR_{++}$ \ and \ $\EE(\zeta) \in \RR_{++}$, \ the condition
 \ $\PP(\zeta > x) = \OO(\PP(\tau > x))$ \ as \ $x \to \infty$ \ does not hold, so Proposition
 \ref{DFK_corr} can not be applied.
Indeed, by Lemma \ref{svosv}, we have \ $\PP(\tau > x) = \oo(\PP(\zeta > x))$ \ as
 \ $x \to \infty$.

We point out that we can not address any result for the tail behavior of the random sum
 \ $\sum_{i=1}^\tau\zeta_i$, \ where \ $\tau$ \ and \ $\zeta$ \ are non-negative regularly varying
 random variables with index \ $\beta \in [1, \infty)$ \ and \ $\alpha \in [1, \infty)$,
 \ respectively such that \ $\alpha < \beta$.
\ If \ $0 \leq \alpha < \beta <1$, \ then Propositions \ref{FGAMS} and \ref{FGAMSRS} can be applied
 as well.

\end{document}